\documentclass[12pt,a4paper]{article}
\usepackage[latin1]{inputenc}
\usepackage{amsmath}
\usepackage{appendix}
\usepackage{multirow}
\usepackage{url}
\usepackage{graphicx}
\usepackage[all]{xy}
\usepackage{amsfonts}
\date{}
\usepackage{amsthm}
\usepackage{amssymb}
\theoremstyle{plain}

\newtheorem{theorem}{Theorem}[section]

\numberwithin{equation}{section}

\begin{document}
\title{Generalising a finite version of Euler's partition identity} 
\author{Darlison Nyirenda} 
\author{Darlison Nyirenda\\ \\
{\small School of Mathematics, University of the Witwatersrand,
Wits 2050, Johannesburg, South Africa.}\\ 
{\small Darlison.Nyirenda@wits.ac.za}  }

\date{}
\maketitle

\begin{abstract}
Recently, George Andrews has given a Glaisher style proof of a finite version of Euler's partition identity. We generalise this result by giving a finite version of Glaisher's partition identity. Both the generating function and bijective proofs are presented.
\bigskip

\noindent 2010 Mathematics Subject Classification: 05A15, 05A17, 05A19, 05A30, 11P81.
\bigskip

\noindent Keywords: partition, identity, bijection, generating function.
\end{abstract}
\section{Introduction}
The so-called Euler's theorem has been widely studied. In the language of integer partitions, the theorem implies that the number of partitions of $n$ into odd parts is equal to the number of partitions of $n$ into distinct parts. This identity, called Euler's partition identity, has been refined (see \cite{eric}).\\
J. W. L. Glaisher gave a bijective proof of the identity (see \cite{andrewsg}). Furthermore, its finite version was given together with bijective proofs (see \cite{corteel}, \cite{bradford} ). We recall this version below.
\begin{theorem}[Euler's theorem-finite version]\label{inde}
The number of partitions of $n$ into odd parts each $\leq 2N$ equals the number of partitions of $n$ into parts each $\leq 2N$ in which the parts $\leq N$ are distinct.
\end{theorem}
For example, if $n = 10$, $N = 3$, then the seven partitions of $n$ into odd parts each $\leq 6$ are\\
$$ 5 + 5,\, 5 + 3 + 1 + 1,\, 5 + 1 + 1 + 1 + 1 + 1,\, 3 + 3 + 3 + 1,$$
$$ 3 + 3 + 1 + 1 + 1 + 1,\, 3 + 1 + 1 + 1 + 1 + 1 + 1 + 1,\, 1 + 1 + 1 + 1 + 1 + 1 + 1 + 1 + 1 + 1.$$
And the seven partitions with parts $\leq 6$ such that parts $\leq 3$ are distinct are\\
$$ 6 + 1,\, 6 + 3 + 1,\, 5 + 5, 5 + 4 + 1,\, 5 + 3 + 2,\, 4 + 4 + 2,\, 4 + 3 + 2 + 1.$$
\noindent However, the bijections for Theorem \ref{inde} given in \cite{bradford} are complicated, and motivated by their complexity, George Andrews gave a much simpler proof that is Glaisher style (see \cite{andrewsg}). \\
It is clear that Euler's partition identity is a specific case of Glaisher's partition identity (see Theorem \ref{glais}) when $s = 2$.
\begin{theorem}[Glaisher's identity, \cite{glaisher}]\label{glais}
The number of partitions of $n$ into parts not divisible by $s$ is equal to the number of partitions of $n$ into parts not repeated more than $s - 1$ times.
\end{theorem}

We are then naturally led to ask as to whether a finite version of Glaisher's partition identity that generalises Theorem \ref{inde} is possible. If so, can we find a bijective proof thereof reminiscent of George Andrews Glaisher style proof?\\

The goal of this paper is to fully address the questions above.\\
Our main result is as follows:
\begin{theorem}\label{indemain}
Let $s$ be a positive integer. The number of partitions of $n$ into parts not divisible by $s$ each $\leq sN$ equals the number of partitions of $n$ into parts each $\leq sN$ in which the parts $\leq N$ occur at most $s - 1$ times.
\end{theorem} 
Observe that Theorem \ref{inde} is obtained by setting $s = 2$.\\
In the subsequent section, we give a generating function proof of the result, and in the section thereafter, a bijective proof that is Glaisher style.
%

\section{First Proof of Theorem \ref{indemain}}
Let $\mathcal{O}_{s,N}(n)$ denote the number of partitions of $n$ in which each part is not divisible by $s$ and $\leq sN$, and $\mathcal{D}_{s,N}(n)$ denote the number of partitions of $n$ in which each part is $\leq sN$ and all parts $\leq N$ occur at most $s - 1$ times. Thus
$$ \sum_{n = 0}^{\infty}\mathcal{O}_{s,N}(n)q^{n} = \prod_{n = 1}^{N}\frac{1}{(1 - q^{sn - 1})(1 - q^{sn - 2})\ldots (1 - q^{sn - s + 1})} $$ and 
$$ \sum_{n = 0}^{\infty}\mathcal{D}_{s,N}(n)q^{n} = \frac{\prod_{n = 1}^{N} (1 + q^{n} + q^{2n} + \ldots + q^{(s - 1)n})}{\prod_{n = 1}^{(s - 1)N}(1 - q^{n + N})}.$$
Observe that
\begin{align*} 
\sum_{n = 0}^{\infty}\mathcal{D}_{s,N}(n)q^{n} & = \frac{\prod_{n = 1}^{N} (1 + q^{n} + q^{2n} + \ldots + q^{(s - 1)n})}{\prod_{n = 1}^{(s - 1)N}(1 - q^{n + N})}\\
                                               & =  \frac{\prod_{n = 1}^{N}(1 - q^{n}) (1 + q^{n} + q^{2n} + \ldots + q^{(s - 1)n})}{\prod_{n = 1}^{N}(1 - q^{n})\prod_{n = 1}^{(s - 1)N}(1 - q^{n + N})}\\
                                               & = \frac{\prod_{n = 1}^{N}(1 - q^{sn})}{\prod_{n = 1}^{sN}(1 - q^{n})}\\
                                               & = \prod_{n = 1}^{N}\frac{1}{(1 - q^{sn - 1})(1 - q^{sn - 2})\ldots (1 - q^{sn - s + 1})}\\
                                               & = \sum_{n = 0}^{\infty}\mathcal{O}_{s,N}(n)q^{n}.
     \end{align*}
\section{Second Proof of Theorem \ref{indemain}}
We give a simple Glaisher style extension of the bijection given by George Andrews \cite{andrewsg}.\\
Consider a partition  $\lambda$ enumerated by $\mathcal{O}_{s,N}(n)$. Each part is of the form $sm - t$ for some $t = 1, 2, \ldots, s - 1$. We can rewrite the partition as
$$\sum_{t = 1}^{s - 1}\sum_{i = 1}^{r_{t}} f_{i,t}(sm_{i} - t) $$
where $f_{i,t}$ is the multiplicity of the part $sm_{i} - t$, $r_{t}$ is the number of parts that are $\equiv -t \pmod{s}$. \\
Note that there exists a unique $\alpha_{i,t}$ such that  $N < (sm_{i} - t)s^{\alpha_{i,t}} \leq sN$. Rather than taking a complete $s$-ary expansion of $f_{i,t}$, we instead do the following:  we find the aforementioned $\alpha_{i,t}$ and use division algorithm to compute $\beta_{i,t}$ and  $e_{i,t}$ from the equation
$$ f_{i,t} = \beta_{i,t}s^{\alpha_{i,t}} + e_{i,t}\,\,\,\text{where}\,\,\,0\leq e_{i,t} \leq s^{\alpha_{i,t}} - 1.$$
Then write the $s$-ary expansion of $e_{i,t}$, i.e.,
$$ e_{i,t} = \sum_{j = 0}^{b_i}a_{j,t}s^{j}\,\,\,\,\text{where}\,\,\,\, 0\leq a_{j,t} \leq s - 1.$$
So $f_{i,t} = \sum_{j = 0}^{b_i}a_{j,t}s^{j} + \beta_{i,t}s^{\alpha_{i,t}}$ and thus
\begin{align*}
\sum_{t = 1}^{s - 1}\sum_{i = 1}^{r_{t}} f_{i,t}(sm_{i} - t) & = \sum_{t = 1}^{s - 1}\sum_{i = 1}^{r_{t}}( \sum_{j = 0}^{b_i}a_{j,t}s^{j} + \beta_{i,t}s^{\alpha_{i,t}})(sm_{i} - t)\\
                               & = \sum_{t = 1}^{s - 1}\sum_{i = 1}^{r_{t}}\sum_{j = 0}^{b_i}a_{j,t}(sm_{i} - t)s^{j}  +  \sum_{t = 1}^{s - 1}\sum_{i = 1}^{r_{t}}\beta_{i,t}(sm_{i} - t)s^{\alpha_{i,t}}\\
                                                              & = \sum_{t = 1}^{s - 1}\sum_{i = 1}^{r_{t}}(a_{0,t}(sm_{i} - t) + a_{1,t}(sm_{i} - t)s + \ldots \\
                              & \hspace{5mm} \,\ldots + a_{b_i,t}(sm_{i} - t)s^{b_i}) + \, \sum_{t = 1}^{s - 1}\sum_{i = 1}^{r_{t}}\beta_{i,t}(sm_{i} - t)s^{\alpha_{i,t}}
\end{align*}

\noindent which is the image of $\lambda$  with parts $(sm_{i} - t)s^{j} \leq N$ that have multiplicity $a_{j,t} \leq s - 1$, and parts $(sm_{i} - t)s^{\alpha_{i,t}} \in [N + 1, sN]$ that have multiplicity $ 0 \leq \beta_{i,t} < \infty$. This image is a partition enumerated by $\mathcal{D}_{s,N}(n)$.\\
The inverse of the bijection is not difficult to contruct. \\
We demontrate the bijection using $n = 177, s = 3, N = 4$ and \\
$\lambda = (11^{6}, 7^{5}, 5^{7},4^{5}, 2^{2}, 1^{17})$.\\
Note that $11.3^{0} = 11.1 \in [5,12]$, $7.3^{0}  = 7.1 \in [5,12]$, $5.3^{0} = 5.1 \in [5,12]$, $4.3^{1} = 4.3 \in [5,12]$, $2.3^{1} = 2.3 \in [5,12]$, and $1.3^{2}  = 1.9 \in [5,12]$.\\
$11^{6}$ is mapped to $11(6.1 + 0) = 6.11$, which is interpreted as $11^{6}$. Similarly, $7^{5}$ to $7^{5}$, and $5^{7}$ to $5^{7}$.\\
$4^{5}$ goes to $4(1.3 + 2)$, using division algorithm $5 = 1.3 + 2$, and now taking the 3-ary expansion of 2; $4(1.3 + 2) = 4(1.3 + 2.3^{0}) = 1.12 + 2.4$, which is $(12, 4^{2})$. Continuing in this manner, we have the mapping  $2^{2} \mapsto 2^{2}$ and $1^{17} \mapsto (9,3^{2},1^{2})$. Thus $$ \lambda \mapsto (12,11^{6},9,7^{5},5^{7},4^{2},3^{2},1^{2}).$$ 


\end{document}